\documentclass[12pt]{article}

\usepackage{amstext    }
\usepackage{amsthm    }
\usepackage{a4}
\usepackage[mathscr]{eucal}
\usepackage{mathrsfs}

\usepackage{amsmath}
\usepackage{amssymb}
\usepackage{amscd}

\newtheorem{theorem}{Theorem}[section]

\newtheorem{proposition}[theorem]{Proposition}
\newtheorem{corollary}[theorem]{Corollary}
\newtheorem{lemma}[theorem]{Lemma}
\newtheorem{remark}[theorem]{Remark}

\newtheorem{example}[theorem]{Example}

\newcommand{\cali}[1]{\mathscr{#1}}

\newcommand{\Aut}{{\rm Aut}}

\newcommand{\id}{{\rm id}}

\newcommand{\Kc}{\cali{K}}

\newcommand{\FS}{{\rm FS}}

\newcommand{\C}{\mathbb{C}}

\newcommand{\R}{\mathbb{R}}
\renewcommand{\P}{\mathbb{P}}

\newcommand{\zwedge}{\stackrel{\circ}{\wedge}}

\title{Comparison of  dynamical degrees for semi-conjugate  meromorphic  maps}

\author{Tien-Cuong Dinh and Vi{\^e}t-Anh Nguy{\^e}n}

\begin{document}

\maketitle

\begin{abstract}
Let $f:X\rightarrow X$ be a dominant meromorphic map on a projective
manifold $X$ which preserves a meromorphic fibration $\pi:X\rightarrow
Y$ of $X$ over a projective manifold $Y$. We establish formulas relating
the dynamical degrees of $f$, the dynamical degrees of $f$ relative to the
fibration and the dynamical degrees of the map $g:Y\rightarrow
Y$ induced by $f$. Applications are given.
\end{abstract}

\noindent
{\bf Classification AMS 2000:} Primary 37F, Secondary  32U40, 32H50.


\section{Introduction} \label{introduction}

Let $(X,\omega_X)$  be  a  compact K{\"a}hler  manifold of dimension
$k$ and let $f:X\rightarrow X$ be a meromorphic map. We assume that $f$ is
{\it dominant}, i.e. the image of $f$ contains an open subset of
$X$. Let $\pi:X\rightarrow Y$ be a dominant meromorphic map from
$X$ onto a compact K{\"a}hler manifold $(Y,\omega_Y)$ of dimension $l\leq
k$. The fibers of $\pi$ define a fibration on $X$ which might be
singular. If $f$ preserves this fibration, i.e. $f$ sends generic
fibers of $\pi$ to fibers of $\pi$, it induces a dominant meromorphic
map $g:Y\rightarrow Y$ such that  $\pi\circ f= g\circ \pi$. 
In that case, we say that $f$ is {\it semi-conjugate} to
$g$. For simplicity, we assume that $\omega_Y$ is so
normalized that $\omega_Y^l$ is a probability measure.

A natural question  is how the dynamical  system defined by $f$ is
similar to the one defined by $g$ when  $f$ is  semi-conjugate  to $g$ as
above. One of the first  steps  towards  understanding  this  question  should be
  to find  out the  relations between some invariants associated to
  $f$ and $g$. In this paper, we will compare their dynamical degrees.

Let $f^n:=f\circ\cdots\circ f$, $n$ times, denote the iterate of order
$n$ of $f$.
The  dynamical degree  $d_p(f)$  of  order $p$ is  the  quantity 
which measures the growth of the  norms of  $(f^n)^*$
acting  on the Hodge cohomology group $H^{p,p}(X,\R)$ when $n$  tends
to infinity. By Poincar{\'e} duality,  it also measures the growth of the  norms of  $(f^n)_*$
acting  on $H^{k-p,k-p}(X,\R)$. If $X$ is a
projective manifold, $d_p(f)$ represents the volume growth of 
 $f^n(V)$ for $p$-dimensional submanifolds $V$  of $X$.  

It was shown by Sibony and the first author in
\cite{DinhSibony1, DinhSibony2} that dynamical degrees are
bi-meromorphic invariants, that is, if $f$ and $g$ are conjugate,
they have the same dynamical degrees. 
Dynamical degrees capture
important dynamical information, in particular, in the computation of
the topological entropy or in the construction of Green currents and
of measures of maximal entropy. We refer  the  reader  to the above
references and to
\cite{DinhSibony5,Gromov,Sibony, Yomdin} for more results on this  matter.

When $f$ preserves a fibration $\pi:X\rightarrow Y$ as above, the
dynamical degree $d_p(f|\pi)$ of order $p$ of $f$ relative to $\pi$ measures the
growth of $(f^n)^*$ acting on the subspace $H_\pi^{l+p,l+p}(X,\R)$
of classes in $H^{l+p,l+p}(X,\R)$ 
which can be supported by a generic fiber of $\pi$. 
It also measures the
growth of $(f^n)_*$ acting on $H_\pi^{k-p,k-p}(X,\R)$ and
represents the volume growth of 
 $f^n(V)$ for $p$-dimensional submanifolds $V$  of a generic fiber of
 $\pi$ when $X$ is projective. Precise definitions and properties will be given in
 Section \ref{section_degree}.
Here is our main result.

\begin{theorem} \label{th_main}
Let  $X$ and $Y$ be projective  manifolds of dimension $k$
and  $l$ respectively with $k\geq l$.
Let $f: X\rightarrow X$, $g:Y\rightarrow Y$ 
and   $\pi:X\rightarrow Y$ be dominant meromorphic maps 
such that $\pi\circ f=g\circ \pi$. 
Then the dynamical degrees $d_p(f)$ of $f$ are related to the dynamical
degrees $d_p(g)$ of $g$ and the relative dynamical degrees $d_p(f|\pi)$
by the formulas
       $$d_p(f)= \max_{\max\{0,p-k+l\}\leq j\leq \min\{p,l\}}d_j(g)d_{p-j}(f|\pi)$$
for $0\leq p\leq k$.
\end{theorem}

Note that the condition $\max\{0,p-k+l\}\leq j\leq \min\{p,l\}$ is
equivalent to $0\leq j\leq l$ and $0\leq p-j\leq k-l$. It guarantees
that $d_j(g)$ and $d_{p-j}(f|\pi)$ are meaningful\footnote{We
  will find later analogous conditions, essentially for the same
  raison but also to avoid expressions which always vanish, e.g. $\omega_Y^{l+1}=0$.}.
We deduce from the above result that $\max d_p(f)\geq \max d_p(g)$.
This gives an affirmative answer to the problem 9.3 in Hasselblatt-Propp \cite{HasselblattPropp}.
When $X$ and $Y$ have the same dimension, generic fibers of $\pi$ are
finite and have the same cardinality. Moreover, $f$ defines bijections
between generic fibers of $\pi$. We deduce from the proof of Theorem \ref{th_main}
the following corollary which generalizes a result in
\cite{DinhSibony1, DinhSibony2}. It was proved by Nakayama-Zhang for
holomorphic maps in \cite{NakayamaZhang}.

 \begin{corollary} \label{cor_degree}
Let  $X$ and $Y$ be compact K{\"a}hler  manifolds of same dimension $k$.
Let $f: X\rightarrow X$, $g:Y\rightarrow Y$ 
and   $\pi:X\rightarrow Y$ be dominant meromorphic maps such that $\pi\circ f=g\circ \pi$.
Then the dynamical degrees of $f$ are equal to the dynamical degrees of
$g$.
\end{corollary}

Recall that by a theorem of Khovanskii \cite{Khovanskii}, Teissier
\cite{Teissier} and Gromov \cite{Gromov1}, the dynamical
degrees of $f$ are log-concave, i.e. $p\mapsto \log d_p(f)$ is
concave. Therefore, there are intergers $p\leq p'$ such that 
$$1=d_0(f)<\cdots<d_p(f)=\cdots=d_{p'}(f)>\cdots>d_k(f).$$
An instructive example with $p\not=p'$ is a map $f(x_1,x_2)=(h(x_1),x_2)$ on a
product $X_1\times X_2$ of projective manifolds. A natural problem is to find dynamically
interesting examples of maps on projective manifolds. For this
purpose, one can try to construct maps with
distinct consecutive dynamical degrees, i.e. with $p=p'$. Somehow, this
condition insures that there is no trivial direction in the associated
dynamical systems. We have the following useful results.

\begin{corollary} \label{cor_distinct_degree}
Let $f,\pi,g$ be as in Theorem \ref{th_main}. If the consecutive
dynamical degrees of $f$ are distinct, then the same property holds
for $g$ and the consecutive dynamical degrees of $f$ relatively to
$\pi$ are also distinct. 
\end{corollary}

The following result is obtained using the Iitaka
fibrations of $X$.

\begin{corollary} \label{cor_kodaira}
Let $X$ be a projective manifold admitting a dominant meromorphic map
with distinct consecutive dynamical degrees. Then the Kodaira
dimension of $X$ is either equal to $0$ or $-\infty$.
\end{corollary} 

Note that the same result  was proved for compact K{\"a}hler surfaces 
by Cantat in \cite{Cantat} and Guedj in
\cite{Guedj}, and for holomorphic maps on compact K{\"a}hler manifolds
by Nakayama and Zhang in \cite{NakayamaZhang,Zhang}. We also refer to
Amerik-Campana \cite{AmerikCampana} and Nakayama-Zhang
\cite{NakayamaZhang, Zhang1} for other invariant fibrations for which Theorem
\ref{th_main} may be applied in order to compute dynamical degrees.

\medskip

\noindent
{\it{\bf Acknowledgment.}}
 The paper was written while  the second  author was visiting  the  Abdus Salam International Centre
 for Theoretical Physics
in Trieste and the Korea Institute for Advanced Study in Seoul. He wishes to express his gratitude to these organizations.

\section{Positive closed currents} \label{section_current}

The proof of our main result uses a delicate calculus on positive
closed currents on compact K\"ahler manifolds\footnote{In this paper, we only consider the strong positivity.}. 
In this section, we prove some useful results which 
can be applied 
to currents of integration on varieties and may have independent
interest. 
The reader will find in Demailly \cite{Demailly} and Voisin 
\cite{Voisin} the basic facts on currents and on K\"ahler geometry.

Let $(X,\omega_X)$ be a
compact K{\"a}hler manifold of dimension $k$. Let $\Kc^p(X)$ denote the
cone of classes of strictly positive closed $(p,p)$-forms in
$H^{p,p}(X,\R)$. This is an open cone which is salient,
i.e. $\overline{\Kc^p(X)}\cap -\overline{\Kc^p(X)}=\{0\}$. 
If $c,c'$ are two classes in $H^{p,p}(X,\R)$, we write $c\leq c'$ and
$c'\geq c$ when $c'-c$ is in $\Kc^p(X)\cup \{0\}$. 

If $T$ is a real
closed $(p,p)$-current, denote by $\{T\}$ its class in
$H^{p,p}(X,\R)$. If moreover $T$ is positive, the {\it mass} of $T$ is defined by
$\|T\|:=\langle T,\omega_X^{k-p}\rangle$. We often use the properties
that $\|T\|$ depends only on the class
of $T$ and $\{T\}\leq
A\{\omega_X^p\}$ for some constant $A>0$ independent of $T$. The
following semi-regularization of currents 
was proved by Sibony and the first author in \cite{DinhSibony1,DinhSibony2}.

\begin{proposition} \label{prop_reg}
Let $T$ be a positive closed $(p,p)$-current on a compact K{\"a}hler
manifold $(X,\omega_X)$. Then there is a sequence of smooth positive
closed $(p,p)$-forms $T_n$ on $X$ which converges weakly to a positive closed
$(p,p)$-current $T'$ such that $T'\geq T$, i.e. $T'-T\geq 0$, $\|T_n\|\leq A\|T\|$ and
$\{T_n\}\leq A\|T\|\{\omega_X^p\}$, where $A>0$ is a constant
independent of $T$. Moreover, if $T$ is smooth on an open set $U$,
then for every compact set $K\subset U$, we have $T_n\geq T$ on $K$
when $n$ is large enough.  
\end{proposition}

Consider now a positive closed $(p,p)$-current $T$ and another
positive closed $(q,q)$-current $S$ on $X$ with $p+q\leq k$. Assume
that $T$ is smooth on  a
dense Zariski open set $U$ of $X$.
Then $T_{|U}\wedge S_{|U}$ is a  well-defined positive closed $(p+q,p+q)$-current
on $U$. The following lemma shows that  $T_{|U}\wedge S_{|U}$ has a
finite mass, 
i.e. $\langle T_{|U}\wedge
S_{|U},\omega_X^{k-p-q}\rangle<+\infty$. Therefore, by Skoda's theorem
\cite{Skoda}, its trivial extension defines a positive closed current on
$X$. We denote by $T\zwedge S$ this current obtained for the maximal
Zariski open set $U$ on which $T$ is smooth. 
Observe that when $S$ has no mass on proper analytic subsets of $X$, the
current obtained in this way does not change if we replace $U$ with another dense
Zariski open set.

\begin{lemma} \label{lemma_wedge_zariski}
Let $T$ and $S$ be as above. Then $T_{|U}\wedge
S_{|U}$ has a finite mass and $T\zwedge S$ is
well-defined. Moreover, we have 
$$\|T\zwedge S\|\leq A\|T\|\|S\|$$
for some constant $A>0$ independent of $T$ and $S$. 
\end{lemma}
\proof
Let $T_n$ and $K$ be as in Proposition \ref{prop_reg}. Since
$\|T_n\wedge S\|$ can be computed cohomologically, we have 
$$\|T_{|K}\wedge S_{|K}\|\leq \liminf_{n\rightarrow\infty} \|T_n\wedge
S\|\leq A\|T\|\|\omega_X^p\wedge S\|=A\|T\|\|S\|.$$
This property holds for every compact subset $K$ of $U$. Therefore, 
$$\|T_{|U}\wedge S_{|U}\|\leq A\|T\|\| S\|.$$
The lemma follows.
\endproof

We will be interested in positive closed currents $T$ on $Y\times \P^m$,
where $(Y,\omega_Y)$ is a compact K{\"a}hler manifold of dimension $l$
and $\P^m$ is the projective space of dimension $m$ endowed with the
standard Fubini-Study form $\omega_\FS$. 
We assume that $\omega_\FS$ is so normalized that $\omega_\FS^m$ is a
probability measure. In practice, we will take $m:=k-l=\dim X-\dim Y$. In order to
simplify the notation, the pull-back of $\omega_Y$ and $\omega_\FS$
to $Y\times \P^m$ under the canonical projections are also denoted by
$\omega_Y$ and $\omega_\FS$. Consider on $Y\times \P^m$ the K{\"a}hler
form $\omega:=\omega_Y+\omega_\FS$. The pull-back of a class $c$ in $H^*(Y,\C)$ or
$H^*(\P^m,\C)$ to $H^*(Y\times \P^m,\C)$ under the canonical projections
is also denoted by $c$.

If $T$ is a positive closed
$(p,p)$-current on $Y\times \P^m$, define for $\max\{0,p-m\}\leq j\leq
\min \{l,p\}$ (or equivalently, for $0\leq j\leq l$ and $0\leq p-j\leq
m$)
$$\alpha_j(T):=\big\langle T,\omega_Y^{l-j}\wedge
\omega_\FS^{m-p+j}\big\rangle.$$
Observe that $\alpha_j(T)$ depends only on the class $\{T\}$ of $T$.
Denote by $\smile$ the cup-product on Hodge cohomology groups.

\begin{proposition} \label{prop_class_prod}
Let $T$ be a positive closed $(p,p)$-current on
  $Y\times\P^m$ as above. Then
$$\{T\}\leq A\sum_{\max\{0,p-m\}\leq j\leq \min \{l,p\}}
\alpha_j(T)\{\omega_Y^j\}\smile \{\omega_\FS^{p-j}\},$$
where $A>0$ is a constant independent of $T$. 
\end{proposition}
\proof
By K{\"u}nneth formula \cite[p.266]{Voisin}, we have
$$H^{\ast}(Y\times \P^m,\C)=H^{\ast}(Y,\C)\otimes H^{\ast}(\P^m,\C).$$ 
Therefore,  there are classes $c_j\in H^{j,j}(Y,\R)$ such that 
\begin{equation*}
\{T\}=\sum_{\max\{0,p-m\}\leq j\leq \min\{l,p\}} c_j\smile \{\omega_\FS^{p-j}\}.
\end{equation*}
Let $S$ be a smooth positive closed $(l-j,l-j)$-form on
$Y$ and $S'$ its canonical pull-back to $Y\times\P^m$. Recall that
$c_j$ denotes also the pull-back of $c_j$ to $Y\times\P^m$.  Since
$\omega_\FS^m$ is a probability measure on $\P^m$, a simple
computation on bidegree gives 
$$c_j\smile\{S\}=  c_j\smile \{S'\}\smile \{\omega_\FS^m\}
= \langle T, S'\wedge \omega_\FS^{m-p+j}\rangle\geq 0.$$

So, $c_j$ belongs to the convex closed cone $\Kc$ of classes $c$ in $H^{j,j}(Y,\R)$ with
$c\smile c'\geq 0$ for $c'\in \Kc^{l-j}(Y)$. Since $\Kc^{l-j}(Y)$ is
open and since $\smile$ is non-degenerate, $\Kc$ is salient,
i.e. $\Kc\cap -\Kc=\{0\}$. The fact
that $\{\omega_Y^{l-j}\}$ is in the interior of $\Kc^{l-j}(Y)$ implies
that $ c_j\smile \{\omega_Y^{l-j}\}=0$ only when $c_j=0$. Moreover, we have
$$\|c_j\|\leq A' c_j\smile \{\omega_Y^{l-j}\}=A'\langle
T,\omega_Y^{l-j}\wedge \omega_\FS^{m-p+j}\rangle=A'\alpha_j(T)$$
for a fixed norm $\|\ \|$ on $H^{j,j}(Y,\R)$ and for some constant
$A'>0$. It follows that 
$$c_j\leq A\alpha_j(T)\{\omega_Y^j\}$$
for some constant $A>0$. The result follows.
\endproof

\begin{proposition} \label{prop_reg_bis}
Let $T$ be a positive closed $(p,p)$-current on $Y\times\P^m$ as
above. Assume that $Y$ is a projective manifold. Then there is a sequence of
smooth positive closed $(p,p)$-forms $T_n$ on $Y\times\P^m$ which
converges weakly to a current $T'\geq T$ such that
$\alpha_j(T_n)\leq A\alpha_j(T)$ for all $j$,
where $A>0$ is a constant independent of $T$. Moreover, if $T$ is
smooth on an open set $U$, then for every compact subset $K$ of $U$
and every $\epsilon>0$, we
have $T_n\geq T-\epsilon \omega^p$ on $K$ when $n$ is large enough.
\end{proposition}
\proof
We first consider the case where $Y=\P^l$ and $\omega_Y$ is the
Fubini-Study form so normalized that $\omega_Y^l$ is a
probability measure. The K{\"u}nneth formula applied to this particular
case says that $T$ is cohomologous to 
$$\sum_{\max\{0,p-m\}\leq j\leq \min \{l,p\}}
\alpha_j(T)\{\omega_Y^j\}\smile \{\omega_\FS^{p-j}\}.$$
Since $Y\times\P^m$ is homogeneous, we can regularize $T$ using the
automorphisms of $Y\times\P^m$ which are close to the identity. 

More precisely, let $\nu_n$ be a sequence of smooth probability measures on
the group of automorphisms $\Aut(Y\times\P^m)$ of $Y\times \P^m$ whose
supports converge to the identity $\id\in\Aut(Y\times\P^m)$. 
Define 
$$T_n:=\int_{\tau\in\Aut(Y\times\P^m)}\tau_*(T)d\nu_n(\tau).$$
Then, $T_n$ are  
smooth positive closed $(p,p)$-forms and converge weakly to
$T$. We also have $\{T_n\}=\{T\}$ and hence 
$\alpha_j(T_n)=\alpha_j(T)$. This gives the first
assertion for $Y=\P^l$. 

For the second assertion, we can prove a stronger property. Let $\Phi$
be a smooth positive $(p,p)$-form on $U$ such that $\Phi\leq T$. We do
not assume that $T$ is smooth nor that $\Phi$ is closed on $U$. Then
$$\Phi_n:=\int_{\tau\in\Aut(Y\times\P^m)}\tau_*(\Phi)d\nu_n(\tau)$$
converge uniformly to $\Phi$ on $K$. Since $\Phi_n\leq T_n$, we have 
$T_n\geq \Phi-\epsilon\omega^p$ on $K$ for $n$ large enough. With our
hypothesis, $T$ is smooth on $U$ and we can replace $\Phi$ with $T$. 

Assume now that $Y$ is a general projective manifold. 
We may find  a finite  family of open holomorphic maps $\Psi_i$,  $1\leq i\leq  s$, from $Y$
onto $\P^l$ such that for every point $y\in Y$ at least one
map $\Psi_i$ is  of maximal rank at $y$.
To do this it suffices to embed $Y$ into a projective space and take a family of central projections.
Let $\Pi_i:\ Y\times \P^m\rightarrow \P^l\times \P^m$  be  defined  by
\begin{equation*}
\Pi_i(y,z):= (\Psi_i(y),z),\quad  (y,z)\in Y\times \P^m.
\end{equation*}
We apply the first case to the currents $T^{(i)}:=(\Pi_i)_{\ast}(T)$.

We construct as above
smooth positive closed $(p,p)$-forms $T^{(i)}_n$ on $\P^l\times
\P^m$ converging to $T^{(i)}$ such that
$\{T^{(i)}_n\}=\{T^{(i)}\}$. Define $T_n:=\sum_i \Pi_i^*(T^{(i)}_n)$. Since the
cohomology classes of $T^{(i)}_n$ are bounded, the classes of $T_n$
are also bounded. Therefore, the masses of $T_n$ are bounded. Up to
extracting a subsequence, we can assume that $\Pi_i^*(T^{(i)}_n)$
converge and hence $T_n$ converge to a
positive closed current $T'$. If $(y,z)$ is a point in $Y\times \P^m$ and $\Psi_i$ has
maximal rank at $y$, then $\Pi_i$ defines a local bi-holomorphic map on
a neighbourhood of $(y,z)$. In this neighbourhood, we have
$$T\leq \Pi_i^*(\Pi_i)_*(T)=\Pi_i^*(T^{(i)})\leq \lim_{n\rightarrow\infty}
\Pi_i^*(T^{(i)}_n)\leq T'.$$
The choice of $\Psi_i$ implies that $T\leq T'$ on $Y\times\P^m$.  
The second assertion of the proposition is a local property. So, it is also easy to check.

It remains to prove the estimate on $\alpha_j(T_n)$. Let
$\widetilde\omega_\FS$ denote the Fubini-Study form of $\P^l$  so
normalized  that 
$\widetilde\omega_\FS^l$ is a probability measure.
Since $\widetilde\omega_\FS$ is strictly positive, there is a constant $A_1>0$ such
that $(\Psi_i)_*\{\omega_Y^{l-j}\}\leq A_1
\{\widetilde\omega_\FS^{l-j}\}$. We also have $(\Psi_i)^*(\widetilde\omega_\FS^{l-j})\leq A_2
\omega_Y^{l-j}$ for some constant $A_2>0$. For simplicity, we will
also denote by $\omega_Y$, $\omega_\FS$ and $\widetilde\omega_\FS$ the
pull-backs of these forms to $Y\times\P^m$ or to $\P^l\times\P^m$. 
In particular, $(\Pi_i)_*(\omega_Y^{l-j}\wedge \omega_\FS^{m-p+j})$
and $(\Psi_i)_{\ast}(\omega_Y^{l-j})\wedge \omega_\FS^{m-p+j}$
represent the same form on $\P^l\times\P^m$.
Since
$T^{(i)}_n$ are smooth and since the following integrals can be computed
cohomologically, we have
\begin{eqnarray*}
\big\langle \Pi_i^{\ast} (T^{(i)}_n),
\omega_Y^{l-j}\wedge \omega_\FS^{m-p+j}\big\rangle
& = &\big\langle T^{(i)}_n,
(\Psi_i)_{\ast}(\omega_Y^{l-j})\wedge \omega_\FS^{m-p+j}\big\rangle\\
& \leq & A_1\big\langle T^{(i)}_n,\widetilde\omega_\FS^{l-j} \wedge \omega_\FS^{m-p+j}\big\rangle\\
& = & A_1\big\langle T^{(i)},\widetilde\omega_\FS^{l-j} \wedge \omega_\FS^{m-p+j}\big\rangle\\
& = & A_1 \big\langle T,
\Psi_i^*(\widetilde\omega_\FS^{l-j}) \wedge \omega_\FS^{m-p+j}\big\rangle\\
&\leq & A_1A_2\big\langle T, \omega_Y^{l-j}\wedge \omega_\FS^{m-p+j}\big\rangle.
\end{eqnarray*}
It follows that $\alpha_j\big(\Pi_i^{\ast} (T^{(i)}_n)\big)\leq A_1A_2
\alpha_j(T)$ and hence $\alpha_j(T_n)\leq A\alpha_j(T)$ for some
constant $A>0$.
\endproof

\section{Dynamical degrees} \label{section_degree}

Let $\pi:(X,\omega_X)\rightarrow (Y,\omega_Y)$ be a dominant
meromorphic map between compact K{\"a}hler manifolds of dimension $k$
and $l$ respectively. The map $\pi$ is holomorphic outside the
indeterminacy set $I_\pi$ which is an analytic
subset of $X$ of codimension at least 2. The closure $\Gamma$ of
its graph over $X\setminus I_\pi$ is an irreducible analytic subset of dimension $k$
of $X\times Y$. If, $\tau_X$ and $\tau_Y$ denote the projections from
$X\times Y$ onto its factors, then $\tau_X$ defines a bi-holomorphic
map between $\Gamma\setminus \tau_X^{-1}(I_\pi)$ and $X\setminus I_\pi$. The
fibers of $\tau_{X|\Gamma}$ over $I_\pi$ have positive dimension. One can
identify $\pi$ with $\tau_Y\circ (\tau_{X|\Gamma})^{-1}$. For $A\subset
X$ and $B\subset Y$,
define $\pi(A):=\tau_Y(\tau_{X|\Gamma})^{-1}(A)$ and 
$\pi^{-1}(B):=\tau_X(\tau_{Y|\Gamma})^{-1}(B)$.

The map $\pi$ induces linear operators on
currents. If $\Phi$ is a smooth $(p,q)$-form on $Y$, then
$\pi^*(\Phi)$ is the $(p,q)$-current defined by
$$\pi^*(\Phi):=(\tau_X)_*(\tau_Y^*(\Phi)\wedge [\Gamma]),$$
where $[\Gamma]$ is the current of integration on $\Gamma$. 
It is not difficult to see that $\pi^*(\Phi)$ is an $L^1$ form smooth
outside $I_\pi$.
If $\Psi$ is a smooth $(p,q)$-form on $X$ with $p,q\geq k-l$, then
$\pi_*(\Psi)$ is the $(p-k+l,q-k+l)$-current defined by
$$\pi_*(\Psi):=(\tau_Y)_*(\tau_X^*(\Psi)\wedge [\Gamma]).$$
If $\Phi$ and $\Psi$ are closed or positive, so are $\pi^*(\Phi)$ and
$\pi_*(\Psi)$. Therefore, $\pi^*$ and $\pi_*$ induce linear operators on the
Hodge cohomology groups of $X$ and $Y$.

In general, the above operators do not extend continuously to positive
closed currents. We will use instead the strict transforms of currents
$\pi^\bullet$ and $\pi_\bullet$ which coincide with $\pi^*$ and
$\pi_*$ on smooth positive closed forms. 
In this paper, we only need these operators in the case where $X$ and
$Y$ have the same dimension $k$. 

Let $U$ be the maximal Zariski open
set in $X\setminus I_\pi$ such that $\pi:U\rightarrow \pi(U)$ is locally invertible.
The complement of $U$ in $X$ is called the {\it critical set} of $\pi$.
If $T$ is a positive closed $(p,p)$-current on $Y$, $(\pi_{|U})^*(T)$
is well-defined and is a positive closed $(p,p)$-current on
$U$. Proposition \ref{prop_reg} allows to show that this current has finite
mass. By Skoda theorem \cite{Skoda}, its trivial extension to $X$
is a positive closed $(p,p)$-current that we denote by
$\pi^\bullet(T)$. 

Let $V$ be the maximal Zariski open set in
$Y\setminus \pi(I_\pi)$ such
that $\pi:\pi^{-1}(V)\rightarrow V$ is a non-ramified covering. The
complement of $V$ in $X$ is called the {\it set of critical values} of
$\pi$. If $S$
is a positive closed $(p,p)$-current on $X$, then $\pi_\bullet(S)$ is
the trivial extension of $(\pi_{|\pi^{-1}(V)})_*(S)$ to $Y$. This is
also a positive closed $(p,p)$-current. We will use the properties that 
$\|\pi^\bullet(T)\|\leq A\|T\|$ and $\|\pi_\bullet(S)\|\leq A\|S\|$
for some constant $A>0$ independent of $T,S$, see \cite{DinhSibony1,
  DinhSibony2} for details. 

Consider now a dominant meromorphic self-map $f:X\rightarrow X$. The
iterate of order $n$ of $f$ is defined by $f^n:=f\circ\cdots\circ f$
($n$ times) on a dense Zariski open set and extends to a
dominant meromorphic map on $X$. Define for $0\leq p\leq k$
$$\lambda_p(f^n):=\|(f^n)^*(\omega_X^p)\|=\big\langle
(f^n)^*(\omega_X^p),\omega_X^{k-p}\big\rangle.$$ 
It is not difficult to see that 
$$\lambda_p(f^n)=\|(f^n)_*(\omega_X^{k-p})\|=\big\langle
(f^n)_*(\omega_X^{k-p}),\omega_X^p\big\rangle.$$ 
It was shown in \cite{DinhSibony1, DinhSibony2} that
$[\lambda_p(f^n)]^{1/n}$ converge to a constant $d_p(f)$ which is 
the {\it dynamical degree of order $p$} of $f$. Note that the main
difficulty here is that in general we do not have
$(f^{n+s})^*=(f^n)^*\circ (f^s)^*$ on currents. 

Let $\|\ \|_{H^{p,p}}$ denote the norm of an operator acting on
$H^{p,p}(X,\R)$ with respect to a fixed norm on that space. Since the
mass of a positive closed current depends only on its cohomology
class, we deduce from the above discussion that
$$A^{-1}\lambda_p(f^n)\leq \|(f^n)^*\|_{H^{p,p}}\leq A\lambda_p(f^n),$$
for some constant $A>0$. It follows that 
$$d_p(f)=\lim_{n\rightarrow\infty}  \|(f^n)^*\|_{H^{p,p}}^{1/n}.$$
Note that we also have $d_p(f^n)=d_p(f)^n$ for $n\geq 1$. 
The last dynamical degree $d_k(f)$ is also called the {\it topological
  degree} of $f$. It is equal to the number of points in a generic
fiber of $f$ and we have $\lambda_k(f^n)\|\omega_X^k\|^{-1}=d_k(f^n)=d_k(f)^n$.

\begin{proposition} \label{prop_current_growth}
Let $T$ be a positive closed $(p,p)$-current and $S$ a positive closed
$(k-p,k-p)$-current on $X$. Then 
$$\|(f^n)^\bullet(T)\|\leq A\|T\|\lambda_p(f^n) \quad \mbox{and} \quad
\|(f^n)_\bullet(S)\|\leq A\|S\|\lambda_p(f^n)$$
for some constant $A>0$ independent of $T$, $S$ and $n$. In particular, we have
$$\limsup_{n\rightarrow\infty} \|(f^n)^\bullet(T)\|^{1/n}\leq d_p(f)
\quad \mbox{and}\quad 
\limsup_{n\rightarrow\infty} \|(f^n)_\bullet(S)\|^{1/n}\leq d_p(f).$$
\end{proposition}
\proof
We show the first inequality. The second one is proved in the same way. 
Let $T_i$ be smooth positive closed forms as in Proposition
\ref{prop_reg}. It follows from the
definition of $(f^n)^\bullet$ that any limit value of $(f^n)^*(T_i)$ is
larger than or equal to $(f^n)^\bullet(T)$. So, it is enough to bound the
mass of $(f^n)^*(T_i)$. Since this mass 
can be computed cohomologically and since 
$\{T_i\}\leq A\|T\|\{\omega_X^p\}$, we obtain that 
$ \|(f^n)^*(T_i)\|\leq A\|T\|\lambda_p(f^n)$ for some constant
$A>0$. This completes the proof.
\endproof

The above proposition can be applied to currents of integration on
submanifolds $V$ of dimension $k-p$ or $p$ of $X$ and gives a upper bound for the
volume growth of the preimage or image of $V$ by $f^n$. 

It was shown in \cite{DinhSibony1, DinhSibony2} that dynamical
degrees are bi-meromorphic invariants, i.e. conjugate maps have the
same dynamical degrees. This property allows to define dynamical
degrees for maps on singular manifolds having a k{\"a}hlerian
desingularization. We will use the same argument in order to define
dynamical degrees relative to an invariant meromorphic fibration. 

Let $\pi:X\rightarrow Y$ be a dominant meromorphic map as in Theorem
\ref{th_main}. It defines a fibration and we assume that $f$ preserves
this fibration. So, $f$ induces a dominant meromorphic map
$g:Y\rightarrow Y$ semi-conjugate to $f$, i.e.  $\pi\circ f=g\circ
\pi$. 
Consider first the case where $\pi$ is a {\bf holomorphic} map. By
Bertini-Sard theorem, the set $Z$  of critical values of $\pi$ is a proper
analytic subset of $Y$. Therefore,
$\pi:X\setminus\pi^{-1}(Z)\rightarrow Y\setminus Z$ defines a regular
holomorphic fibration. Its fibers form a continuous family of 
smooth submanifolds of dimension $k-l$ of $X$. 

Let $P_f$ and $P_g$ denote the union of the critical set and the set
of critical values of $f$ and $g$ respectively. They contain the
indeterminacy sets of $f$ and of $g$. 
A fiber $L_y:=\pi^{-1}(y)$
with $y\in Y\setminus Z$ is called {\it generic} if for every
$n\geq 0$
\begin{enumerate}
\item[(a)] $g^n(y)$ and $g^{-n}(y)$ do not intersect $P_g$;
\item[(b)] For every point $b$ in $g^n(y)\cup g^{-n}(y)$, no
  component of $L_b$ is contained in $P_f$.
\end{enumerate}
Denote by $\Sigma$ the set of $y$ such that $L_y$ is generic. Observe
that $Y\setminus \Sigma$ is contained in a finite
or countable union of proper analytic subsets of $Y$. 
So, $\Sigma$ is connected. We also have
$g(\Sigma)=g^{-1}(\Sigma)=\Sigma$. We will use the following lemma for
$\nu=\omega_Y^l$ and for $\nu=[d_l(g)]^{-n}(g^n)^*(\omega_Y^l)$.

\begin{lemma} \label{lemma_deg_pi}
Let $L_y$ be a generic fiber as above. Let $\nu$ be a probability
measure on $Y$ which has no mass on proper analytic subsets of $Y$. Then, for $0\leq p\leq k-l$
and for $n\geq 0$, the 6 positive closed currents 
$$d_l(g)^{-n}(f^n)^\bullet(\omega_X^p\wedge [L_y]),\quad
(f^n)^*(\omega_X^p)\zwedge [L_y],\quad 
(f^n)^*(\omega_X^p)\zwedge \pi^*(\nu)
$$ 
and
$$(f^n)_\bullet(\omega_X^{k-l-p}\wedge
[L_y]),\quad 
d_l(g)^{-n}(f^n)_*(\omega_X^{k-l-p})\zwedge
[L_y],\quad
d_l(g)^{-l}(f^n)_*(\omega_X^{k-l-p}) \zwedge \pi^*(\nu)$$ 
have the same mass. In particular, their mass does
not depend on $y\in\Sigma$.
\end{lemma}
\proof
For $y\in\Sigma$, define
$$\varphi(y):=d_l(g)^{-n}\big\|(f^n)^\bullet(\omega_X^p\wedge [L_y])\big\|
\quad \mbox{and} \quad 
\psi(y):=\big\|(f^n)_\bullet(\omega_X^{k-l-p}\wedge
[L_y])\big\|.$$
It is not difficult to see that these functions are continuous on $\Sigma$. 
We have
\begin{eqnarray*}
\varphi(y) & = & d_l(g)^{-n}\big\langle (f^n)^\bullet(\omega_X^p\wedge [L_y]),\omega_X^{k-l-p}\big\rangle \\
& = &  d_l(g)^{-n}\big\langle\omega_X^p, [L_y]\zwedge
(f^n)_*(\omega_X^{k-l-p})\big\rangle.
\end{eqnarray*}
It follows that 
$$\varphi= d_l(g)^{-n}\pi_*\big(\omega_X^p\wedge (f^n)_*(\omega_X^{k-l-p})\big)$$
in the sense of currents on $Y$. Therefore, $\varphi$ defines a closed $0$-current
on $Y$ and it should be constant on $\Sigma$.

We also deduce from the above computation that 
$$\varphi(y)= d_l(g)^{-n}\big\|(f^n)_*(\omega_X^{k-l-p})\zwedge [L_y]\big\|.$$ 
Since $\nu$ has no mass on $Y\setminus\Sigma$, we
obtain
$$\varphi =\int \varphi(y) d\nu =  d_l(g)^{-n}\big\|(f^n)_*(\omega_X^{k-l-p})\zwedge \pi^*(\nu)\big\|.$$

In the same way, we prove that $\psi$ is constant on $\Sigma$ and
$$\psi = \big\|(f^n)^*(\omega_X^p)\zwedge
[L_y]\big\|=\big\|(f^n)^*(\omega_X^p)\zwedge \pi^*(\nu)\big\|.$$
It remains to check that $\varphi=\psi$. Using that $\psi$ is
constant and $\#g^{-n}(y)=d_l(g)^n$, we have
$$\varphi=d_l(g)^{-n}(f^n)^\bullet(\omega_X^p\wedge [L_y]) =
d_l(g)^{-n}\sum_{b\in g^{-n}(y)}(f^n)^*(\omega_X^p)
  \zwedge [L_b]=\psi.$$
This completes the proof.
\endproof

Define $\lambda_p(f^n|\pi)$ the mass of the currents in Lemma
\ref{lemma_deg_pi}. We have in particular
$$\lambda_p(f^n|\pi)=\big\|(f^n)^*(\omega_X^p)\wedge \pi^*(\omega_Y^l)\big\|.$$

\begin{proposition} \label{prop_def_deg_pi}
The sequence $\lambda_p(f^n|\pi)^{1/n}$ converges
  to a constant $d_p(f|\pi)$. Let $T$ be a positive closed
  $(p+l,p+l)$-current and $S$ a positive closed $(k-p,k-p)$-current on
  $X$ which are 
supported on a generic fiber $L_y$. Then 
$$\|(f^n)^\bullet(T)\|\leq A_y\|T\|d_l(g)^n\lambda_p(f^n|\pi) \quad \mbox{and} \quad
\|(f^n)_\bullet(S)\|\leq A_y\|S\|\lambda_p(f^n|\pi)$$
for some constant $A_y>0$ independent of $T$ and $S$. In particular, we have
$$\limsup_{n\rightarrow\infty} \|(f^n)^\bullet(T)\|^{1/n}\leq d_l(g)d_p(f|\pi)
\quad \mbox{and}\quad 
\limsup_{n\rightarrow\infty} \|(f^n)_\bullet(S)\|^{1/n}\leq d_p(f|\pi).$$
\end{proposition}
\proof
Fix a generic fiber $L_y$ with $y\in\Sigma$. 
We will show that 
$$\lambda_p(f^{n+m}|\pi)\leq A_y
\lambda_p(f^n|\pi)\lambda_p(f^m|\pi)$$ 
for some constant $A_y>0$ and for
all $n,m\geq 0$. This will imply the first assertion because the
sequence $A_y\lambda_p(f^n|\pi)$ is sub-multiplicative.

Since $L_y$ is a compact K{\"a}hler manifold, we
can apply Proposition \ref{prop_reg} to $L_y$. Let $b$ be a point in
$\Sigma$ such that $g^m(b)=y$. Define
$R:=(f^m)_\bullet(\omega_X^{k-l-p}\wedge [L_b])$. 
This is a positive closed $(k-p,k-p)$-current on $X$ which is also a
$(k-l-p,k-l-p)$-current on $L_y$. 
By Lemma
\ref{lemma_deg_pi}, we have $\|R\|= \lambda_p(f^m|\pi)$. Therefore,
there are smooth positive closed $(k-l-p,k-l-p)$-forms $\Theta_i$ on
$L_y$ which converge to a current $\Theta\geq R$. Moreover, we have
$\{\Theta_i\}\leq A_y \lambda_p(f^m|\pi)\{\omega_{X|L_y}^{k-l-p}\}$ for some constant
$A_y>0$, where the inequality is considered in $H^*(L_y,\R)$. 

Let $h$ denote the restriction of $f^n$ to $L_y$. It defines
a meromorphic map from $L_y$ to $L_{g^n(y)}$. Since the mass of a
positive closed current can be computed cohomologically, we obtain
\begin{eqnarray*}
\lambda_p(f^{n+m}|\pi) &= & \|(f^n)_\bullet(R)\|\leq \liminf_{i\rightarrow\infty}
\|h_*(\Theta_i)\| \leq A_y
\lambda_p(f^m|\pi)\|h_*(\omega_{X|L_y}^{k-l-p})\|\\
& = & A_y \lambda_p(f^m|\pi)\|(f^n)_\bullet(\omega_X^{k-l-p}\wedge
[L_y])\| = A_y \lambda_p(f^m|\pi)\lambda_p(f^n|\pi).
\end{eqnarray*}
This implies the first assertion in the proposition.
The rest is proved in the same way using the semi-regularization
result for $T$ and $S$ on
$L_y$, see also Proposition \ref{prop_current_growth}.
\endproof

We call $d_p(f|\pi)$ the {\it dynamical degree of order $p$} of $f$
{\it relatively} to $\pi$. 
The convergence in Proposition \ref{prop_def_deg_pi}
implies that $d_p(f^n|\pi)=d_p(f|\pi)^n$.

\begin{remark}\rm
Our choice of $\Sigma$ simplifies the calculus on currents but several
properties above still hold for some $y$ out of $\Sigma$. For example,
if $y$ is a fixed point of $g$ which is not a critical value of $\pi$
and if no component of $L_y$ is contained in the critical set of
$f$, then $d_p(f|\pi)=d_p(f_{|L_y})$. The proof is left to the
reader. 
\end{remark}

The next result shows that the relative dynamical degrees are 
bi-meromorphic invariants. Consider a bi-meromorphic map
$\tau:(\widetilde X,\omega_{\widetilde X})\rightarrow (X,\omega_X)$ between compact
K{\"a}hler manifolds. Define $\widetilde \pi:=\pi\circ\tau$ and
$\widetilde f:=\tau^{-1}\circ f\circ \tau$. Then, $\widetilde f$ is a dominant
meromorphic map conjugate to $f$ and $\widetilde\pi\circ \widetilde
f=g\circ\widetilde\pi$. 

\begin{proposition} \label{prop_deg_rel_bim}
Assume that $\widetilde\pi$ is holomorphic. Then
$$d_p(f|\pi)=d_p(\widetilde f|\widetilde\pi)$$
for $0\leq p\leq k-l$.
\end{proposition}
\proof
Since $\tau$ is bi-meromorphic,
$\tau_\bullet\widetilde\pi^*(\omega_Y^l)= \pi^*(\omega_Y^l)$
and $\widetilde
f^n=\tau^{-1}\circ f^n\circ \tau$, we have
\begin{eqnarray*}
\lambda_p(\widetilde f^n|\widetilde\pi) & = & 
\big\langle (\widetilde
f^n)^*(\omega_{\widetilde X}^p)\wedge \widetilde\pi^*(\omega_Y^l),
\omega_{\widetilde X}^{k-l-p}\big\rangle \\
& = &  \big\langle \tau^\bullet(f^n)^\bullet \tau_*(\omega_{\widetilde
  X}^p)\wedge \omega_{\widetilde X}^{k-l-p}, \widetilde\pi^*(\omega_Y^l)
\big\rangle  \\
& = & \big\langle (f^n)^\bullet \tau_*(\omega_{\widetilde
  X}^p)\zwedge \tau_*(\omega_{\widetilde X}^{k-l-p}),
\pi^*(\omega_Y^l)\big\rangle. 
\end{eqnarray*}
Using the semi-regularization result for $\tau_*(\omega_{\widetilde X}^{k-l-p})$, we
deduce that 
$$\lambda_p(\widetilde f^n|\widetilde\pi)\leq A \big\langle (f^n)^\bullet \tau_*(\omega_{\widetilde
  X}^p)\wedge \omega_X^{k-l-p},
\pi^*(\omega_Y^l)\big\rangle$$
for some constant $A>0$. Then, using a semi-regularization of
$\tau_*(\omega_{\widetilde X}^p)$, we obtain
$$\lambda_p(\widetilde f^n|\widetilde\pi)\leq A' \big\langle (f^n)^*(\omega_X^p)\wedge \omega_X^{k-l-p},
\pi^*(\omega_Y^l)\big\rangle =A' \lambda_p(f^n|\pi)$$
for some constant $A'>0$. It follows that $d_p(\widetilde
f|\widetilde\pi)\leq d_p(f|\pi)$. The converse inequality is proved in
the same way.
\endproof

The last proposition allows to define relative dynamical degrees in
the general case. Assume now that $f$ preserves a {\bf meromorphic}
fibration $\pi:X\rightarrow Y$, i.e. $\pi\circ f=g\circ\pi$ as in
Theorem \ref{th_main}. Let $\Gamma$ denote the closure of the graph of
$\pi$ in $X\times Y$. Then $\Gamma$ is an irreducible analytic set of
dimension $k$ which is bi-meromorphic to $X$. Let $\sigma:\widetilde
X\rightarrow \Gamma$ be a desingularization of $\Gamma$ which can be
constructed using a blow-up
along the singularities. By Blanchard's
theorem \cite{Blanchard}, $\widetilde X$ is a compact K{\"a}hler
manifold. Then, $\tau:=\tau_X\circ \sigma$ is a bi-meromorphic map from
$\widetilde X$ to $X$. Define also $\widetilde\pi:=\tau_Y\circ \sigma$ and
$\widetilde f:=\tau^{-1}\circ f\circ\tau$. The map $\widetilde\pi$ is
holomorphic and $\widetilde\pi\circ \widetilde
f=g\circ\widetilde\pi$. Define the {\it dynamical degree of order $p$} of $f$
{\it relatively} to $\pi$ by
$$d_p(f|\pi):=d_p(\widetilde f|\widetilde\pi).$$
Proposition \ref{prop_deg_rel_bim} implies that the definition does not depend on the
choice of $\sigma$. The following result is a consequence of a theorem by
Khovanskii, Teissier and Gromov.

\begin{proposition} The function $p\mapsto \log d_p(f|\pi)$ is concave
  for $0\leq p\leq k-l$. In particular, $d_p(f|\pi)\geq 1$ for $0\leq
  p\leq k-l$.  
\end{proposition}
\proof
We can assume that $\pi$ is holomorphic.
We have to show that 
$$d_{p-1}(f|\pi)d_{p+1}(f|\pi)\leq d_p(f|\pi)^2.$$ 
For this purpose, it is enough to check that 
$$\lambda_{p-1}(f^n|\pi)\lambda_{p+1}(f^n|\pi)\leq \lambda_p(f^n|\pi)^2.$$ 
Observe that for non-critical values $y$ of $\pi$, the fibers $L_y$ are not necessarily connected but
they contain the same number $s$ of components. The family of these
components is connected since $X$ is connected. It defines a covering
of degree $s$ over the set of non-critical values of $\pi$. Let
$\Sigma'$ denote the set of the components $L_y'$ of $L_y$ with $y\in\Sigma$.
We can prove as in Lemma \ref{lemma_deg_pi}  
that the function 
$L_y'\mapsto \|(f^n)^*(\omega_X^p)\zwedge [L_y']\|$
is constant on $\Sigma'$. 
Therefore, it is
equal to $s^{-1}\|(f^n)^*(\omega_X^p)\zwedge [L_y]\|$ and
then to $s^{-1}\lambda_p(f^n|\pi)$. 

Let $h$ be the restriction of $f^n$ to $L_1:=L_y'$ and define
$L_2:=h(L_1)$. Let $\Gamma$ denote the graph of $h$ in $L_1\times L_2$
and $\tau:\widehat\Gamma\rightarrow\Gamma$ a desingularization of
$\Gamma$ using some blow-up along the singularities. By Blanchard's
theorem \cite{Blanchard}, $\widetilde\Gamma$ is a compact K{\"a}hler
manifold. Denote by $\tau_1:\widehat\Gamma\rightarrow L_1$ and
 $\tau_2:\widehat\Gamma\rightarrow L_2$ the canonical projections. We
 have $h=\tau_2\circ\tau_1^{-1}$. Define
 $\omega_1:=\tau_1^*(\omega_X)$ and $\omega_2:=\tau_2^*(\omega_X)$. 
We deduce from the above
 discussion that
$$s^{-1}\lambda_p(f^n|\pi)=\|(h^n)^*(\omega_X^p)\|=\int_{\widehat\Gamma}
\omega_1^{k-l-p}\wedge\omega_2^p.$$
If $\gamma_p$ denotes the last integral, Gromov proved in
\cite{Gromov1} that $p\mapsto\log \gamma_p$ is concave,
i.e. $\gamma_{p-1}\gamma_{p+1}\leq\gamma_p^2$, when 
$\omega_1$ and $\omega_2$ are  K{\"a}hler forms. By continuity, this
still holds in our
case where these forms are only smooth positive and closed. Hence,
$p\mapsto \log d_p(f|\pi)$ is concave. 

In order to deduce the second assertion of the proposition, it is
enough to show that $d_0(f|\pi)=1$ and $d_{k-l}(f|\pi)\geq 1$. 
For $y$ generic, we have 
$$\lambda_0(f^n|\pi)=d_l(g)^{-n}\|(f^n)^\bullet[L_y]\|=d_l(g)^{-n}\sum_{b\in
  g^{-n}(y)}\|[L_b]\|.$$
Hence, $\lambda_0(f^n|\pi)$ is independent of $n$
since $\#g^{-n}(y)=d_l(g)^n$ and the mass of $[L_b]$, with
$b\in\Sigma$, is independent of $b$. It follows that $d_0(f|\pi)=1$.

We also have for $y$ generic and $b\in g^{-n}(y)$
$$\lambda_{k-l}(f^n|\pi)=\|(f^n)_\bullet[L_b]\|\geq \|[L_y]\|.$$
So, the sequence $\lambda_{k-l}(f^n|\pi)$ is bounded from below by a
positive constant. Therefore, $d_{k-l}(f|\pi)\geq 1$. This completes
the proof of the lemma. Note that we can show that $d_{k-l}(f|\pi)$ is
the number of points in a generic fiber of the restriction of $f$ to $L_y$.
\endproof

Consider now some examples, see also \cite{AmerikCampana, NakayamaZhang,
  Zhang, Zhang1}.

\begin{example} \rm
Let $X=Y\times Z$ be the product of two compact K\"ahler manifolds and
$\pi:X\rightarrow Y$ the canonical projection. Consider $f(y,z):=(g(y),h(z))$
where $g:Y\rightarrow Y$ and $h:Z\rightarrow Z$ are dominant
meromorphic maps. So, $f$ is semi-conjugate to $g$. The relative dynamical degree
$d_p(f|\pi)$ is equal to $d_p(h)$. We easily deduce from the
definition of dynamical degrees that
$$d_p(f)=\max_{\max\{0,p-k+l\}\leq j\leq \min\{p,l\}}
d_j(g)d_{p-j}(h).$$
There are more interesting examples of maps on the product $Y\times
Z$. Let $F$ be a compact K\"ahler manifold. Assume that $F$ is also the parameter
space of a meromorphic family of meromorphic self-maps of $Z$. 
Let $\tau:Y\rightarrow F$ be a meromorphic map.
Then $f(y,z):=(g(y),\tau(y)(z))$ is a meromorphic self-map of $Y\times Z$ which
preserves the fibration $\pi$. The example is also interesting 
when $\tau(y)$ is holomorphic for generic $y$ or when a Zariski open
set $G$ of $F$ is a Lie group and $\tau$ is a morphism from $G$ to the
group of bi-meromorphic maps of $Z$.
\end{example}

\begin{example} \rm
Let $g:Y\rightarrow Y$ be a dominant meromorphic map on a compact
K\"ahler manifold $Y$. It induces a meromorphic self-map $f$ on the
projectivization $X:=\P T_Y$ of the holomorphic tangent bundle of
$Y$. The map $f$ preserves the fibration associated to the canonical projection from $X$ onto $Y$
and is semi-conjugate to $g$. This example and some applications were
considered in \cite{Dinh}.
\end{example}

\section{Proofs of the results} \label{section_proof}

We first prove Theorem \ref{th_main}. Since the dynamical degrees are
bi-meromorphic invariants, we can assume that $\pi$ is a holomorphic
map. Since $X$ is projective, we can construct a dominant meromorphic map
$v:X\rightarrow \P^{k-l}$. Indeed, it is enough to embed $X$ in a
projective space and choose a generic central projection on $\P^{k-l}$. 
Replacing $X$ with a desingularization of the graph of $v$ allows to assume that $v$ is holomorphic. 
Consider the
holomorphic map $\Pi:X\rightarrow Y\times\P^{k-l}$ defined by
$$\Pi(x):=(\pi(x),v(x)).$$  
Since the chosen central projection is generic, the intersection of a
generic fiber of $\pi$ and a generic fiber of $v$ is
finite. Therefore, $\Pi$ is dominant. 

Our proof is based on a delicate calculus on currents. If
$X=Y\times\P^{k-l}$ and $\pi$ is the canonical projection onto $Y$,
the proof is simpler and the properties obtained in Section
\ref{section_current} can be directly applied. A rough idea is to
reduce the general case to the particular case using the map $\Pi$. In
other words, we use the fact that $f$ is, in some sense, ``semi-conjugate'' to
the multi-valued map $\Pi\circ f\circ\Pi^{-1}$ which is defined on
$Y\times \P^{k-l}$.

Let $\omega_\FS$ denote the Fubini-Study form on $\P^{k-l}$. For simplicity, the canonical pull-back of
$\omega_Y$ and $\omega_\FS$ to $Y\times
\P^{k-l}$ are still denoted by   $\omega_Y$ and $\omega_\FS$. 
In particular, $\Pi^*(\omega_Y)$ and $\pi^*(\omega_Y)$ represent the
same form on $X$. 
We consider on $Y\times \P^{k-l}$ the K{\"a}hler form
$\omega:=\omega_Y+\omega_\FS$. Our calculus will involve the
quantities $a_{q,p}(n)$ defined
for $n\geq 0$,
$0\leq q\leq k-l$ and $q\leq p\leq l+q$ by
\begin{eqnarray*}
a_{q,p}(n)  & := &   \|\Pi_*(f^n)^* \Pi^*(\omega^p)\wedge \omega_Y^{l-p+q}\|=\big\langle \Pi_*(f^n)^*
\Pi^*(\omega^p),   \omega_Y^{l-p+q}\wedge \omega^{k-l-q}\big\rangle \\
& = & \big\langle (f^n)^*
\Pi^*(\omega^p),   \Pi^*(\omega_Y^{l-p+q}\wedge
\omega^{k-l-q})\big\rangle \\
& = & \big\langle (f^n)^*
\Pi^*(\omega^p)\wedge   \pi^*(\omega_Y^{l-p+q}),
\Pi^*(\omega^{k-l-q})\big\rangle.
\end{eqnarray*}
Observe that 
$$a_{q,p}(n)  \geq \alpha_{p-q}\big(\Pi_*(f^n)^*\Pi^*(\omega^p)\big),$$
where $\alpha_{p-q}(\cdot)$ are defined in Section \ref{section_current}.

\begin{lemma} \label{lemma_a_pp}
There is a constant $A>0$ independent of $p,n$ such that
$$A^{-1} \lambda_p (f^n|\pi)\leq a_{p,p}(n)\leq  A\lambda_p
(f^n|\pi)$$ 
In particular, $[a_{p,p}(n)]^{1/n}$ converge to $d_p(f|\pi)$.
\end{lemma}
\proof
Since the pull-back of a smooth form under $\Pi$ is smooth, we have
\begin{eqnarray*}
a_{p,p}(n) & = &  \big\langle (f^n)^*
\Pi^*(\omega^p)\wedge \pi^*(\omega_Y^{l}),
\Pi^*(\omega^{k-l-p})\big\rangle \\
& \leq &  A\big\langle
(f^n)^*(\omega_X^p)\wedge\pi^*(\omega_Y^l),\omega_X^{k-l-p}\big\rangle =
A\lambda_p (f^n|\pi)
\end{eqnarray*}
for some constant $A>0$. This gives the second inequality in the lemma.

Define $T:=\Pi_*(\omega_X^p)$. Since
$\Pi^\bullet(T)\geq \omega_X^p$, we have
$$\lambda_p(f|\pi)=\|(f^n)^*(\omega_X^p)\wedge \pi^*(\omega_Y^l)\|
\leq \|(f^n)^\bullet\Pi^\bullet(T)\wedge \pi^*(\omega_Y^l)\|.$$
We apply Proposition \ref{prop_reg} to the current
$T$ on $Y\times\P^{k-l}$ which is an $L^1$ form smooth on a Zariski open set. Let $T_i$ be as in 
that proposition 
with $\{T_i\}\leq A\{\omega^p\}$ for some constant $A>0$. 
If $S:=\Pi_*(\omega_X^{k-l-p})$, we have
$\Pi^\bullet(S)\geq\omega_X^{k-l-p}$ and hence
\begin{eqnarray*}
\lambda_p(f|\pi) & \leq &  \liminf_{i\rightarrow\infty }
\|(f^n)^*\Pi^*(T_i)\wedge \pi^*(\omega_Y^l)\| \leq
A\|(f^n)^*\Pi^*(\omega^p)\wedge \pi^*(\omega_Y^l)\| \\
& = & A\big\langle (f^n)^*\Pi^*(\omega^p)\wedge \pi^*(\omega_Y^l),
\omega_X^{k-l-p}\big\rangle\leq 
A\big\|(f^n)^*\Pi^*(\omega^p)\wedge \pi^*(\omega_Y^l)\zwedge
\Pi^\bullet(S)\big\|.
\end{eqnarray*}

Now, we apply again Proposition \ref{prop_reg}, in particular its last
assertion, to the current $S$ which is an $L^1$ form smooth on a Zariski open
set. If $S_i$ are
smooth forms satisfying that proposition, the latter expression is
bounded from above by
$$\liminf_{i\rightarrow\infty}  \big\langle (f^n)^*\Pi^*(\omega^p)\wedge \pi^*(\omega_Y^l),
\Pi^*(S_i)\big\rangle \lesssim \big\langle (f^n)^*\Pi^*(\omega^p)\wedge \pi^*(\omega_Y^l),
\Pi^*(\omega^{k-l-p})\big\rangle.$$
The last integral is equal to $a_{p,p}(n)$.
The first inequality in the lemma follows.
\endproof

Define for $0\leq p\leq k$
$$b_p(n):=\sum_{\max\{0,p-l\}\leq q\leq \min\{p,k-l\}} a_{q,p}(n).$$
We have the following lemma.

\begin{lemma} \label{lemma_b_p}
The sequence $b_p(n)^{1/n}$ converges to $d_p(f)$.
\end{lemma}
\proof
Since $\Pi^*(\omega^p)$, $\pi^*(\omega_Y^{l-p+q})$
and $\Pi^*(\omega^{k-l-q})$ are smooth on $X$, we have 
$$a_{q,p}(n)= \big\langle (f^n)^*
\Pi^*(\omega^p)\wedge   \pi^*(\omega_Y^{l-p+q}),
\Pi^*(\omega^{k-l-q})\big\rangle \leq
A \|(f^n)^*(\omega_X^p)\|=A\lambda_p(f^n)$$
for some constant $A>0$. We deduce that $\limsup b_p(n)^{1/n}\leq
d_p(f)$. 

It remains to check that $\liminf b_p(n)^{1/n}\geq d_p(f)$. For this
purpose, we only need to show that $\lambda_p(f^n)\leq Ab_p(n)$ for
some constant $A>0$. 
Define $T:=\Pi_*(f^n)^* \Pi^*(\omega^p)$. We prove that 
$\lambda_p(f^n)\lesssim \|T\|\lesssim b_p(n)$ which will imply the result.

Define $S:=\Pi_*(\omega_X^p)$. We have $\Pi^\bullet(S)\geq
\omega_X^p$. Therefore,
$$\lambda_p(f^n)=\big\langle (f^n)^*(\omega_X^p),\omega_X^{k-p}\big\rangle \leq
\big\langle (f^n)^\bullet\Pi^\bullet(S),\omega_X^{k-p}\big\rangle.$$
Using a semi-regularization of $S$, we deduce that
$$\lambda_p(f^n)\lesssim \big\langle (f^n)^*\Pi^*(\omega^p),\omega_X^{k-p}\big\rangle.$$
Define $R:=\Pi_*(\omega_X^{k-p})$. We also have $\Pi^\bullet(R)\geq
\omega_X^{k-p}$. We obtain as above using a semi-regularization of $R$ that 
\begin{eqnarray*}
\lambda_p(f^n) & \lesssim &  \big\|
(f^n)^*\Pi^*(\omega^p)\zwedge\Pi^\bullet(R)\big\|\lesssim \big\|
(f^n)^*\Pi^*(\omega^p)\wedge\Pi^*(\omega^{k-p})\big\|\\
& = & \big\langle \Pi_*(f^n)^*\Pi^*(\omega^p),
\omega^{k-p}\big\rangle=\|T\|.
\end{eqnarray*}

Now, since $\omega_Y^{l+1}=0$ and $\omega_\FS^{k-l+1}=0$, we have
\begin{eqnarray*}
\|T\| & = & \big\langle T, (\omega_Y+\omega_\FS)^{k-p}\big\rangle \lesssim 
\sum_{\max\{0,p-l\}\leq q \leq \min\{p,k-l\}} \big\langle T,\omega_Y^{l-p+q}\wedge
\omega_\FS^{k-l-q} \big\rangle \\
& \leq &   \sum_{\max\{0,p-l\}\leq q \leq \min\{p,k-l\}} a_{p,q}(n) = b_p(n).
\end{eqnarray*}
This completes the proof of the lemma.
\endproof

For every $n\geq 0$ and $0\leq p\leq l$ define  
$$c_p(n) := \lambda_p(g^n)=\|(g^n)^{*}(\omega_Y^p)\|=\big\langle
(g^n)^{*}(\omega_Y^p), \omega_Y^{l-p}\big\rangle.$$
We have the following lemma.

\begin{lemma} \label{lemma_Pi_wedge}
There is a constant $A>0$ such that
$$\big\langle \Pi_*
(f^n)^*\Pi^*(\omega_Y^{p-q}\wedge\omega^q),\omega_Y^{l-p+p_0}\wedge
\omega^{k-l-p_0}\big\rangle \leq A a_{p_0,q}(n)c_{p-q}(n)$$
for $0\leq p_0\leq k-l$, $p_0\leq p\leq l+p_0$, $p_0\leq q\leq p$ and $n\geq
0$.
Moreover, the above integral vanishes when $q<p_0$.  
\end{lemma}
\proof
Observe that by definition of $\Pi_*$
\begin{eqnarray*}
\Pi_*(f^n)^*\Pi^*(\omega_Y^{p-q}\wedge\omega^q) 
& = & \Pi_*\big[(f^n)^*\Pi^*(\omega_Y^{p-q})\zwedge
(f^n)^*\Pi^*(\omega^q)\big] \\
& \leq &  \Pi_* (f^n)^*\Pi^*(\omega_Y^{p-q})\zwedge\Pi_*
(f^n)^*\Pi^*(\omega^q).
\end{eqnarray*}
Hence, the left hand side of the inequality in the lemma is smaller
than or
equal to
$$\big\langle \Pi_*
(f^n)^*\Pi^*(\omega_Y^{p-q})\zwedge\Pi_*
(f^n)^*\Pi^*(\omega^q),\omega_Y^{l-p+p_0}\wedge
\omega^{k-l-p_0}\big\rangle.$$
Define $T:=\Pi_*(f^n)^*\Pi^*(\omega_Y^{p-q})\wedge \omega_Y^{l-p+p_0}$
and $S:=\Pi_*(f^n)^*\Pi^*(\omega^q)\wedge\omega^{k-l-p_0}$. Note that
$T$ and $S$ are of bidegree $(l-q+p_0,l-q+p_0)$ and
$(k-l+q-p_0,k-l+q-p_0)$ respectively. 
The quantity considered above is equal to the mass of the measure
$T\zwedge S$. 

We first show that
$\alpha_j(T)=0$ when $j<l-q+p_0$ and $\alpha_{l-q+p_0}(T)\leq A
c_{p-q}(n)$ for some constant $A>0$. Since
$\pi\circ f^n=g^n\circ\pi$, we have
$$T= \Pi_*(f^n)^*\pi^*(\omega_Y^{p-q})\wedge \omega_Y^{l-p+p_0}
=\Pi_* \pi^\bullet(g^n)^*(\omega_Y^{p-q})\wedge \omega_Y^{l-p+p_0}.$$
Hence,
\begin{eqnarray*}
\alpha_j(T) & = &  \big\langle \Pi_*
\pi^\bullet(g^n)^*(\omega_Y^{p-q})\wedge \omega_Y^{l-p+p_0},
\omega_Y^{l-j}\wedge \omega_\FS^{k-2l+q-p_0+j}\big\rangle \\
& = &  \big\langle 
\pi^\bullet(g^n)^*(\omega_Y^{p-q})\wedge \pi^*(\omega_Y^{l-p+p_0}),
\pi^*(\omega_Y^{l-j})\wedge \Pi^*(\omega_\FS^{k-2l+q-p_0+j})\big\rangle
\\
& = &  \big\langle 
\pi^\bullet\big[(g^n)^*(\omega_Y^{p-q})\wedge \omega_Y^{2l-p+p_0-j}\big],
\Pi^*(\omega_\FS^{k-2l+q-p_0+j})\big\rangle. 
\end{eqnarray*}
When $j<l-q+p_0$, the form in the brackets has bidegree $\geq (l+1,l+1)$
and should vanish because $\dim Y=l$. Therefore,  $\alpha_j(T)=0$ in that
case. When $j=l-q+p_0$,
this form defines a positive measure of mass
$\lambda_{p-q}(g^n)$. Its cohomology class is equal to $\lambda_{p-q}(g^n)\{\omega_Y^l\}$.
Therefore, using a semi-regularization as above, we obtain
$$\alpha_{l-q+p_0}(T)  \lesssim   \lambda_{p-q}(g^n) \big\langle
\pi^*(\omega_Y^l),\Pi^*(\omega_\FS^{k-l})\big\rangle
\leq Ac_{p-q}(n)$$
for some constant $A>0$.

We deduce from Proposition \ref{prop_class_prod} that $\{T\}\lesssim
c_{p-q}(n)\{\omega_Y^{l-q+p_0}\}$. Using the semi-regularization
proposition \ref{prop_reg_bis} for $T$, we obtain
$$\{T\zwedge S\} \lesssim c_{p-q}(n)\|\omega_Y^{l-q+p_0}\wedge S\|=
c_{p-q}(n)a_{p_0,q}(n).$$
This completes the proof of the first assertion in the lemma.
For the second one, it is enough to observe 
that when $q<p_0$, we have
$\alpha_j(T)=0$ for every $j$ and hence $T=0$. 
\endproof

The following lemma is crucial in our proof.

\begin{lemma}\label{lemma_key}
 There exists a constant $A>0$ such that for all $0\leq p_0\leq k-l$, $p_0\leq p\leq l+p_0$ and all $n,r\geq 1$
$$ a_{p_0,p}(nr)
\leq  A^r \sum \prod_{s=1}^r a_{p_{s-1},p_s}(n)c_{p-p_s}(n),$$
where the  sum is taken over $p_s$ with 
$p_{s-1}\leq p_s\leq p$ and $p_{s-1}\leq k-l$ for $s=1,\ldots,r$.
\end{lemma}
\proof
We proceed  by induction on $r$. Clearly, the lemma is true for
$r=1$. 
Suppose  the lemma true for $r$, we need to prove  it for $r+1$. In
what follows, the constants $A_i$ depend only on the geometry
of $X$ and $Y$. 

Define $T^{(r)}:=\Pi_* (f^{nr})^*\Pi^*(\omega^p)$. This is a
positive closed $L^1$ form, smooth on a dense Zariski open set.
Observe that $\Pi^\bullet\Pi_*\geq \id$ on positive closed currents
having no mass on proper analytic subsets of $X$. Therefore,
$$T^{(r+1)}\leq
\Pi_*(f^n)^\bullet\Pi^\bullet\Pi_*(f^{nr})^\bullet\Pi^*(\omega^p) = 
\Pi_* (f^n)^\bullet\Pi^\bullet(T^{(r)}).$$

On the other hand, by Proposition \ref{prop_reg_bis}, we can find a
sequence of smooth positive closed $(p,p)$-forms $T_i^{(r)}$
converging  weakly to a positive closed current $\widetilde T^{(r)}\geq T^{(r)}$
such that
$$\alpha_{p-q}(T_i^{(r)})\leq A_1\alpha_{p-q}(T^{(r)}) \leq  A_1 a_{q,p}(nr)$$
for $\max\{0,p-l\}\leq q\leq \min\{p,k-l\}$ and $A_1>0$ a constant.
By Proposition \ref{prop_class_prod}, there is  a constant $A_2>0$
such that
$$\{T_i^{(r)}\}\leq A_2 \sum_{\max\{0,p-l\}\leq q \leq \min\{p,k-l\}}
a_{q,p}(nr)\{\omega_Y^{p-q}\}\smile \{\omega_\FS^q\}.$$

We deduce from the above discussion and Lemma \ref{lemma_Pi_wedge} that 
\begin{eqnarray*}
\lefteqn{a_{p_0,p}(n(r+1))   =   \big\langle T^{(r+1)},\omega_Y^{l-p+p_0}\wedge
\omega^{k-l-p_0}\big\rangle} \\
& \leq & \liminf_{i\rightarrow\infty} \big\langle \Pi_*
(f^n)^*\Pi^*(T^{(r)}_i),\omega_Y^{l-p+p_0}\wedge
\omega^{k-l-p_0}\big\rangle \\
& \leq & A_2 \sum_{\max\{0,p-l\}\leq q\leq \min\{p,k-l\}} a_{q,p}(nr)\big\langle \Pi_*
(f^n)^*\Pi^*(\omega_Y^{p-q}\wedge\omega_\FS^q),\omega_Y^{l-p+p_0}\wedge
\omega^{k-l-p_0}\big\rangle \\
& \leq & A_3 \sum_{p_0\leq q\leq\min\{p,k-l\}} a_{q,p}(nr) a_{p_0,q}(n)c_{p-q}(n)
\end{eqnarray*}
for some constant $A_3>0$. Consequently,
the induction hypothesis implies the result.
\endproof

Theorem \ref{th_main} is a consequence of the next two propositions.

\begin{proposition} We have
 $$d_p(f)\geq  \max_{\max\{0,p-k+l\}\leq j\leq
   \min\{p,l\}}d_j(g)d_{p-j}(f|\pi)$$
for $0\leq p\leq k$.
\end{proposition}
\proof
Since $\Pi^*(\omega_Y^j\wedge \omega^{p-j})$ is a smooth form, we have for some constant $A>0$
$$\big\|(f^n)^*\Pi^*(\omega_Y^j\wedge \omega^{p-j})\big\|\leq A\lambda_p(f^n).$$
So, by definition of dynamical degrees and Lemma \ref{lemma_a_pp}, 
it is enough to bound  $\|(f^n)^*\Pi^*(\omega_Y^j\wedge \omega^{p-j})\|$ from below by a constant times 
$\lambda_j(g^n) a_{p-j,p-j}(n)$. 

Fix a constant $A>0$ large enough. Using the identity $\pi\circ
f^n=g^n\circ\pi$ and that 
$\Pi^*(\omega_Y^{l-j}\wedge 
\omega^{k-l-p+j})$ is smooth, we obtain
\begin{eqnarray*}
\lefteqn{A\|(f^n)^*\Pi^*(\omega_Y^j\wedge \omega^{p-j})\|} \\ 
& \geq & \big\langle  (f^n)^*\Pi^*(\omega_Y^j\wedge \omega^{p-j}), \Pi^*(\omega_Y^{l-j}\wedge 
\omega^{k-l-p+j})\big\rangle\\
& = & \big\langle  (f^n)^*\pi^*(\omega_Y^j)\zwedge (f^n)^*\Pi^*(\omega^{p-j}), \pi^*(\omega_Y^{l-j})\wedge 
\Pi^*(\omega^{k-l-p+j})\big\rangle\\
& = & \big\|  (f^n)^*\pi^*(\omega_Y^j)\wedge \pi^*(\omega_Y^{l-j}) \zwedge (f^n)^*\Pi^*(\omega^{p-j})\wedge 
\Pi^*(\omega^{k-l-p+j})\big\|\\
& = & \big\|  \pi^*[(g^n)^*(\omega_Y^j)\wedge \omega_Y^{l-j}]\zwedge (f^n)^*\Pi^*(\omega^{p-j})\wedge 
\Pi^*(\omega^{k-l-p+j})\big\|.
\end{eqnarray*}

Observe that $(g^n)^*(\omega_Y^j)\wedge \omega_Y^{l-j}$ is a positive measure of mass $\lambda_j(g^n)$. As in 
Lemma \ref{lemma_deg_pi}, we show that the last expression  
is equal to $\lambda_j(g^n)$ times the mass of the restriction of $(f^n)^*\Pi^*(\omega^{p-j})\wedge 
\Pi^*(\omega^{k-l-p+j})$ to a generic fiber $L_y$ of $\pi$. Therefore, it is also equal to
$$\lambda_j(g^n)\big\langle  \pi^*(\omega_Y^l), (f^n)^*\Pi^*(\omega^{p-j})\wedge 
\Pi^*(\omega^{k-l-p+j})\big\rangle=\lambda_j(g^n)a_{p-j,p-j}(n).$$
This completes the proof.
\endproof

\begin{proposition} We have
 $$d_p(f)\leq  \max_{\max\{0,p-k+l\}\leq j\leq
   \min\{p,l\}}d_j(g)d_{p-j}(f|\pi)$$
for $0\leq p\leq k$.
\end{proposition}
\proof
For every $0\leq  p\leq k$ and $n\geq 0$ let
$$\mu_p(n):=\max_{\max\{0,p-k+l\}\leq j\leq
   \min\{p,l\}} c_j(n)a_{p-j,p-j}(n).$$
Observe that for $r>p$,  in Lemma \ref{lemma_key},
there are at most $p$ indices $s$  such that  $p_{s-1}<p_s$.
Moreover, the sum in that lemma contains at most
$(k+1)^r$ terms and the sum in the definition of $b_p(n)$ contains
at most $p+1$ terms. 
We infer the following estimate
$$b_p(rn)\leq \Big[(p+1)(k+1)^rA^r b_0(n)\cdots b_p(n)\prod_{j=0}^l c_j(n)\Big]
\mu_p(n)^r.$$
We deduce that
$$[b_p(rn)]^{1/rn}\leq (p+1)^{1/nr}(k+1)^{1/n}A^{1/n} \big[b_0(n)^{1/n}\cdots b_p(n)^{1/n}\big]^{1/r}
\Big[\prod_{j=0}^l c_j(n)^{1/n}\Big]^{1/r}
\mu_p(n)^{1/n}.$$
Letting $n$ tend to infinity, we obtain using Lemma \ref{lemma_b_p} that
$$d_p(f)  \leq   \big[d_0(f)\ldots d_p(f)\big]^{1/r} 
\Big[\prod_{j=0}^l d_j(g)\Big]^{1/r}\liminf_{n\to\infty} \mu_p(n)^{1/n} .$$
Now, letting $r\rightarrow\infty$, the first two factors in the right
hand  side tend to 1. Therefore,  using Lemma \ref{lemma_a_pp}, 
we obtain 
$$d_p(f)\leq  \liminf_{n\to\infty} \mu_p(n)^{1/n} =\max_{\max\{0,p-k+l\}\leq j\leq
   \min\{p,l\}}d_j(g)d_{p-j}(f|\pi).$$
This completes the proof.
\endproof

\noindent
{\bf Proof of Corollary \ref{cor_degree}.}
When $X$ and $Y$ are projective, the corollary is a direct consequence of Theorem
\ref{th_main}. We only used the projectivity in Proposition
\ref{prop_reg_bis} applied to $m:=k-l$ and for the existence of $v:X\rightarrow\P^{k-l}$. This 
is superfluous when $X$ and $Y$ have the same dimension, i.e. $k=l$.
\hfill $\square$

\bigskip

\noindent
{\bf Proof of Corollary \ref{cor_distinct_degree}.}
Let $j$ and $p$ be such that $d_j(g)=\max_q d_q(g)$ and
$d_{p-j}(f|\pi)=\max_q d_q(f|\pi)$. 
We have $0\leq j\leq l$ and $0\leq p-j\leq k-l$. 
By Theorem \ref{th_main}, $d_p(f)$
is the maximal dynamical degree of $f$ and
$d_p(f)=d_j(g)d_{p-j}(f|\pi)$. 
We have
$d_{p-1}(f)<d_p(f)<d_{p+1}(f)$. 
Theorem \ref{th_main} implies that 
$$d_{j-1}(g)<d_j(g)<d_{j+1}(g)
\quad \mbox{and}\quad
d_{p-j-1}(f|\pi)<d_{p-j}(f|\pi)<d_{p-j+1}(f|\pi).$$
The log-concavity of $d_q(g)$ and $d_q(f|\pi)$ implies the
result. Note that when $j=0,l$ or $p-j=0,k-l$, in the above
inequalities, one has to remove the expressions
which are not meaningful.
\hfill $\square$

\bigskip

In the rest of the paper, we prove Corollary \ref{cor_kodaira}. 
Let $K_X$ denote the canonical lines bundle of $X$. Let $H^0(X,K_X^n)$
denote the space of holomorphic sections of $K_X^n$ and 
$H^0(X,K_X^n)^*$ its dual space. Assume that $H^0(X,K_X^n)$ has a
positive dimension. If $x$ is a generic point in $X$, the family $H_x$ of
sections which vanish at $x$ is a hyperplane of  $H^0(X,K_X^n)$
passing through 0. So,
the correspondence $x\mapsto H_x$ defines 
a meromorphic map
$$\pi_n:X\rightarrow \P H^0(X,K_X^n)^*$$
from $X$ to the projectivization of $H^0(X,K_X^n)^*$ which is called an
{\it Iitaka fibration} of $X$. Let $Y_n$ denote the image of $X$ by
$\pi_n$. The {\it Kodaira dimension} of $X$ is
$\kappa_X:=\max_{n\geq 1}\dim Y_n$. When $H^0(X,K_X^n)=0$ for
every $n\geq 1$, the Kodaira dimension of $X$ is $-\infty$. 
We have the following result.

\begin{theorem}[\cite{NakayamaZhang, Ueno}] \label{thm_Iitaka_fix}
Let $f:X\rightarrow X$ be a dominant meromorphic map. Assume that
$\kappa_X\geq 1$. Then $f$ preserves the Iitaka fibration
$\pi_n:X\rightarrow Y_n$. Moreover, the map $g:Y_n\rightarrow Y_n$
induced by $f$ is periodic, i.e. $g^N=\id$ for some integer $N\geq
1$. 
\end{theorem}

\noindent
{\bf Proof of Corollary \ref{cor_kodaira}.}
Assume in order to get a contradiction that $\kappa_X\geq 1$. Let $n\geq 1$ be such that $l:=\dim Y_n\geq
1$. Replacing $f$ with an iterate, we can assume that $g=\id$. A
priori, $Y_n$ may be singular, but we can use a blow-up and assume
that $Y_n$ is smooth. We have $d_j(g)=1$ for $0\leq j\leq l$.
This contradicts
Corollary \ref{cor_distinct_degree}. 
Note that in order to prove that $d_j(g)=1$, instead of Theorem \ref{thm_Iitaka_fix}, 
it is enough to use the weaker result that $g$ is induced
by a linear endomorphism of $\P H^0(X,K_X^n)^*$. 
\hfill $\square$

\noindent
T.-C. Dinh, UPMC Univ Paris 06, UMR 7586, Institut de
Math{\'e}matiques de Jussieu, 4 place Jussieu, F-75005 Paris,
France.\\ 
{\tt  dinh@math.jussieu.fr}, {\tt http://www.math.jussieu.fr/$\sim$dinh}

\medskip

\noindent
V.-A.  Nguy{\^e}n,
Vietnamese Academy  of Science  and  Technology,
Institute of Mathematics,
Department  of Analysis,
18  Hoang Quoc  Viet  Road, Cau Giay  District,
10307 Hanoi, Vietnam.
{\tt nvanh@math.ac.vn}

\noindent
{\sc Current  address:}
 School of Mathematics,
Korea Institute  for Advanced Study,
 207-43 Cheongryangni-2dong,
Dongdaemun-gu,
Seoul 130-722, Korea. {\tt vietanh@kias.re.kr}

\

\end{document}